\documentclass[12pt]{article}
\usepackage{amssymb,amsfonts,amsmath,amsthm}
\usepackage{comment}
\usepackage[normalem]{ulem}
\usepackage[colorinlistoftodos,prependcaption,textsize=tiny]{todonotes}

\usepackage{color}
\usepackage{graphicx}

\newcommand{\beq}{\begin{eqnarray*}}
\newcommand{\feq}{\end{eqnarray*}}
\newcommand{\beqn}{\begin{eqnarray}}
\newcommand{\feqn}{\end{eqnarray}}

\newcommand{\la}{{\lambda}}

\makeatletter \@addtoreset{theorem}{section}\makeatother

\newtheorem*{theorem*}{Theorem}

\font\header=cmssdc10 at 20pt

\begin{document}

{\header Implosion of a pure death process}

\vskip1cm

Luiz Renato Fontes,
Instituto de Matem\'atica e Estatística,
Universidade de S\~ao Paulo,
Rua do Mat\~ao 1010,
05508-090 São Paulo SP 
Brasil
\newline email: lrfontes@usp.br

Rinaldo B. Schinazi, Department of Mathematics, University of Colorado, Colorado Springs, CO 80933-7150, USA;\newline e-mail: 
Rinaldo.Schinazi@uccs.edu

\vskip1cm

{\bf Abstract.} We study a pure death process. At each discrete time every individual dies or not independently of each other with a constant probability.
We give examples showing that in a certain limit extinction happens along a path where one and only one individual is lost at a time. We also exhibit an example 
for which such a path goes from infinity to 0 in a finite time. This is what we call a process implosion.

\bigskip

{\bf 1. The model}

\bigskip

For species on a death spiral there are many examples for which the species ends when a single last individual dies. There was a last known passenger pigeon,
a last Aurochs, a last Tasmanian tiger and so on... A new word ("endling") has even been proposed to designate the last individual of a species. In this paper we propose
a probability model to compute the chance that a species ends up with a single individual. Next we describe our model.

Consider the following discrete time pure death process $(D_t)$. Let $c\in (0,1)$ be a fixed parameter. 
At every time $t\geq 0$ we think of $D_t$ as the number of individuals alive at time $t$. Every individual alive at time $t$
dies with probability $c$ at time $t+1$, independently of each other. In other words,
if $D_t=x$ where $x$ is a positive integer then 
$$D_{t+1}=x-bin_t(x,c),$$
where $bin_t(x,c)$ is a binomial random variable with parameters $x$ and $c$. The different binomial random variables involved
in the construction of the process are independent. Note that $0$ is an absorbing state for this process. This model is a particular case of a model
that goes back to at least Neuts (1994), see also Ben-Ari, Roitershtein and Schinazi (2017). 

\bigskip

{\bf 2. The extinction time}

\bigskip

With a population starting with $n$ individuals let $\tau_n$ be the time for the process to hit 0. Observe that
$$P(\tau_n\leq t)=\left(1-(1-c)^t\right)^n.$$
From this equality it is easy to prove the following limit in probability,
$$\lim_{n\to\infty}\frac{\tau_n}{d_n}=1,$$
where 
$$d_n=-\frac{\ln n}{\ln(1-c)}.$$
Note that $\tau_n$ is actually the maximum of $n$ i.i.d. geometric random variables with parameter $c$. There is a large literature on the limiting behavior of
such a sequence of  random variables, see for instance Eisenberg (2008).

\bigskip

{\bf 3. A special path to extinction}

\bigskip

Starting with $k$ individuals at time  $t=0$ the process may stay at $k$ or drop to any $j$ for $0\leq j\leq k-1$ at time $t=1$.
Let $A_k$ be the event that when the process leaves state $k$ it drops from $k$ to $k-1$. Conditioning on the number of times the process stays at $k$ before dropping we get 
$$P(A_k)=\sum_{j=1}^\infty \left((1-c)^k\right)^{j-1}k(1-c)^{k-1}c.$$
Hence,
$$P(A_k)=\frac{k(1-c)^{k-1}c}{1-(1-c)^k}.\leqno (1)$$

We are interested in the following special path to extinction. Let $B_n$ be the event that starting at $n$ the process gets to 0 only with drops of exactly one. That is,
$$B_n=\bigcap_{k=1}^n A_k.$$
By the independence of the $A_k$ we get
\begin{align*}
P(B_n)&=\prod_{k=1}^n P(A_k)\\
&=\prod_{k=1}^n \frac{k(1-c)^{k-1}c}{1-(1-c)^k}.
\end{align*}

For fixed $c$ it is easy to see that $\lim_{n\to\infty} P(B_n)=0$. 
We will next give examples showing that $\lim_{n\to\infty} P(B_n)>0$ provided
$c$ is a sequence converging to 0 fast enough. The following inequality will be useful.
Since $(1-c)^k\geq 1-kc$ for all $c$ in $[0,1]$ and all $k\geq 1$, then
$$P(B_n)\geq \prod_{k=1}^n (1-c)^{k-1}.$$

\medskip

$\bullet$ Assume that $c$ depends on the initial state $n$ of the process. Then,
\begin{align*}
P(B_n)\geq &\prod_{k=1}^n (1-c_n)^{k-1}\\
& =(1-c_n)^{\frac{n(n-1)}{2}}\\
\end{align*}
If the sequence $(n^2 c_n)$ has a finite limit then $\lim_{n\to\infty} P(B_n)>0$.
If $(n^2 c_n)$ converges to 0 then $\lim_{n\to\infty} P(B_n)=1$.

\medskip

We now turn to our second example.

\medskip

$\bullet$ Assume that $c$ depends on the current state of the process. Then,
$$P(B_n)\geq \prod_{k=1}^n (1-c_k)^{k-1}.$$
If the series $\sum_{k=1}^\infty kc_k<+\infty$ then  $\lim_{n\to\infty} P(B_n)>0$.

\medskip

$\bullet$ Assume that $c$ depends on both the initial state $n$ and the current state $k$. Consider the particular case
$$c_{k,n}=\frac{k^{\alpha}}{n^{\beta}},$$
for $1\leq k\leq n$ where $\alpha>0$ and $\beta>0$ are parameters. Note that
\begin{align*}
P(B_n)\geq & \prod_{k=1}^n \left (1-\frac{n^{\alpha}}{n^{\beta}}\right)^{k-1}\\
=& \left (1-\frac{1}{n^{\beta-\alpha}}\right)^{\frac{n(n-1)}{2}}
\end{align*}
Hence, for $\beta-\alpha>2$ we get that $\lim_{n\to\infty} P(B_n)=1$.

\bigskip

{\bf 4. First passage times}

\bigskip

For $k\geq 1$ let $T_k$ be the time for the process to drop from $k$ to $k-1$. This need not happen so there is a positive probability that $T_k=+\infty$.
Let $s>0$ and
$$g_k(s)=E(e^{sT_k};T_k<\infty).$$
For every $j\geq 1$,
$$P(T_k=j)=\left((1-c)^k\right)^{j-1}k(1-c)^{k-1}c.$$
Hence,
\begin{align*}
g_k(s)=&\frac{kce^s(1-c)^{k-1}}{1-e^s(1-c)^k}\\
=&\frac{kc(1-c)^{k-1}}{e^{-s}-(1-c)^k}
\end{align*}

Next we use the expression above in two examples.

$\bullet$ Let $(a_n)$ be a positive sequence going to infinity and $(c_n)$ be a sequence in $(0,1)$ going to $0$ so that
$$\lim_{n\to\infty}a_nc_n=\la>0.$$
For $n\geq 1$, consider the process for $c=c_n$. Then $g_k$ depends on $n$ and we write $g_{k,n}$ for $g_k$. For $k\geq 1$,
$$g_{k,n}(\frac{s}{a_n})=\frac{kc_n(1-c_n)^{k-1}}{e^{-s/a_n}-(1-c_n)^k}.$$
For fixed $k\geq 1$,
$$\lim_{n\to\infty}ka_nc_n(1-c_n)^{k-1}=k\la\mbox { and }\lim_{n\to\infty}a_n\left(e^{-s/a_n}-(1-c_n)^k\right)=-s+k\la.$$
Therefore,
$$\lim_{n\to\infty}g_{k,n}(\frac{s}{a_n})=\frac{k\la}{k\la-s}.$$
This shows that $\frac{T_{k,n}}{a_n}$ on $\{T_{k,n}<+\infty\}$ converges in distribution to an exponential distribution with parameter $k\la.$

\medskip

We now turn to our second example.

\medskip

$\bullet$ Let $c_{k,n}$ be in $(0,1)$ for every $k\geq 1$ and $n\geq 1$. Assume also that for any fixed $k$,
$$\lim_{n\to\infty}\frac{c_{k,n}}{k^{\alpha}/n^{\beta}}=1,$$
for  $\beta>0$ and $\alpha>0$.
Let $a_n=n^\beta$.
$$g_{k,n}(\frac{s}{a_n})=\frac{kc_{k,n}(1-c_{k,n})^{k-1}}{e^{-s/a_n}-(1-c_{k,n})^k}.$$
Using that
$$\lim_{n\to\infty}ka_nc_{k,n}n(1-c_{k,n})^{k-1}=k^{\alpha+1}\mbox { and }\lim_{n\to\infty}\left(a_ne^{-s/a_n}-a_n(1-c_{k,n})^k\right)=-s+k^{\alpha+1},$$
$\frac{T_{k,n}}{a_n}$ on $\{T_{k,n}<+\infty\}$ converges in distribution to an exponential distribution with parameter $k^{\alpha+1}.$
Note that by (1) the probability of the event $\{T_{k,n}<+\infty\}$ approaches 1
as $n$ goes to infinity.

For this limiting process the expected time to go from $+\infty$ to $0$ is
$$\sum_{k=1}^\infty \frac{1}{k^{\alpha+1}}.$$
Therefore, this process goes from infinity to $0$ in a finite time! In this sense it can be said that the process implodes. A similar behavior is observed for a different model by Menshikov and Petritis (2012).

\vskip1cm

{\bf References}

\bigskip

I. Ben-Ari, A. Roitershtein and R.B. Schinazi (2017) A random walk with catastrophes.

B.~Eisenberg (2008). On the expectation of the maximum of IID geometric random variables.
Statist. Probab. Lett. 78, 135-143.

M. Menshikov and D. Petritis (2012) Explosion, implosion, and moments of passage times for continuous-time Markov chains: a semi-martingale approach.

https://arxiv.org/abs/1202.0952

M.~F.~Neuts (1994). An interesting random walk on the non-negative integers.
J. Appl. Probab. 31, 48-58.

\end{document}